\documentclass[a4paper,11pt]{article}
\usepackage{amsmath, amssymb, amsthm}
\usepackage{graphicx} 

\usepackage[latin1]{inputenc}

\begin{document}

\def\Z{{\mathbb Z}}


\title{\bf{Basic questions on Artin-Tits groups}}

\author{\textsc{Eddy Godelle\footnotemark[1] and Luis Paris\footnote{Both authors are partially supported by the {\it Agence Nationale de la Recherche} ({\it projet Théorie de Garside}, ANR-08-BLAN-0269-03).}}}

\date{\today}

\maketitle

\begin{abstract}
\noindent
This paper is a short survey on four basic questions on Artin-Tits groups: the torsion, the center, the word problem, and the cohomology ($K(\pi,1)$ problem). It is also an opportunity to prove three new results concerning these questions: (1) if all free of infinity Artin-Tits groups are torsion free, then all Artin-Tits groups will be torsion free; (2) If all free of infinity irreducible non-spherical type Artin-Tits groups have a trivial center then all irreducible non-spherical type Artin-Tits groups will have a trivial center; (3) if all free of infinity Artin-Tits groups have solutions to the word problem, then all Artin-Tits groups will have solutions to the word problem. Recall that an Artin-Tits group is free of infinity if its Coxeter graph has no edge labeled by $\infty$.
\end{abstract}



\section{Introduction}

Let $S$ be a finite set. A {\it Coxeter matrix} over $S$ is a square matrix $M=(m_{s,t})_{s,t\in S}$ indexed by the elements of $S$ such that $m_{s,s}=1$ for all $s \in S$, and $m_{s,t} = m_{t,s} \in \{2,3,4, \dots, \infty\}$ for all $s,t \in S$, $s \neq t$. A Coxeter matrix is usually represented by its {\it Coxeter graph}, $\Gamma$. This is a labeled graph whose set of vertices is $S$, such that two vertices $s,t \in S$ are joined by an edge if $m_{s,t} \ge 3$, and such that this edge is labeled by $m_{s,t}$ if $m_{s,t} \ge 4$.

\bigskip\noindent
The {\it Coxeter system} associated with $\Gamma$ is the pair $(W,S)$ where $W=W_\Gamma$ is the group presented by 
\[
W = \left\langle S \left| \begin{array}{cl}
s^2=1& \text{for } s \in S\\
(st)^{m_{s,t}}=1&\text{for } s,t \in S, s \neq t \text{ and } m_{s,t} \neq \infty
\end{array} \right. \right\rangle\,.
\]
The group $W$ is called {\it Coxeter group} of type $\Gamma$.

\bigskip\noindent
If $a,b$ are two letters and $m$ is an integer greater or equal to $2$, we set $\Pi(a,b : m) = (ab)^{\frac{m}{2}}$ if $m$ is even, and $\Pi(a,b : m) = (ab)^{\frac{m-1}{2}} a$ if $m$ is odd. Let $\Sigma= \{\sigma_s ; s \in S\}$ be an abstract set in one-to-one correspondence with $S$. The {\it Artin-Tits system} associated with $\Gamma$ is the pair $(A,\Sigma)$, where $A=A_\Gamma$ is the group presented by 
\[
A=\langle \Sigma\ |\ \Pi(\sigma_s, \sigma_t : m_{s,t}) = \Pi(\sigma_t, \sigma_s : m_{s,t}) \text{ for } s,t \in S, s \neq t \text{ and } m_{s,t} \neq \infty \rangle\,.
\]
The group $A$ is called {\it Artin-Tits group} of type $\Gamma$.

\bigskip\noindent
Coxeter groups were introduced by Tits in his manuscript \cite{Tits1}, which was used by Bourbaky as a basis for writing his seminal book ``Lie groups and Lie algebras, Chapters 4, 5, and 6'' \cite{Bourb1}. They are involved in several areas of mathematics such as group theory, Lie theory, or hyperbolic geometry, and are in some sense fairly well understood. We recommend \cite{Humph1}, \cite{Davis1} and \cite{Bourb1} for a detailed study on these groups.

\bigskip\noindent
Artin-Tits groups were also introduced by Tits \cite{Tits2}, as extensions of Coxeter groups. There is no general result on these groups, and the theory consists on the study of more or less extended families. In particular, the following basic questions are still open. 
\begin{itemize}
\item[{\bf (1)}]
Do Artin-Tits groups have torsion?
\item[{\bf (2)}]
How is the center of Artin-Tits groups?
\item[{\bf (3)}]
Do Artin-Tits groups have solutions to the word problem? 
\item[{\bf (4)}]
What can we say about the cohomology of Artin-Tits groups?
\end{itemize}
Nevertheless, there are fairly precise conjectures on these four questions, and the history of the theory of Artin-Tits groups is intimately related to them. They will be presented within they historical context in Section 2.

\bigskip\noindent
A Coxeter graph $\Gamma$ is called {\it free of infinity} if $m_{s,t} \neq \infty$ for all $s,t \in S$, $s \neq t$. In Section 3 we prove the following principle on the conjectures presented in Section 2 (except for the one on the cohomology, which has been already proved in \cite{EllSko1} and \cite{GodPar1}). 

\bigskip\noindent
{\bf Principle.}
{\it If a property is true for all free of infinity Artin-Tits groups, then it will be true for all Artin-Tits groups.}

\bigskip\noindent
A careful reader may point out to the authors that this principle is false for the property ``to be free of infinity'', but we are sure that he is clever enough to understand that this principle is a general idea and not a theorem, and it may have exceptions. 

\bigskip\noindent
Finally, note that the four questions mentioned above are not the only questions discussed in the theory of Artin-Tits groups. There are many others, more or less important, such as the conjugacy problem, the automaticity, the linearity, or the orderability. 


\section{Conjectures}

The statement of the conjecture on the torsion is quite simple:

\bigskip\noindent
{\bf Conjecture A.}
{\it Artin-Tits groups are torsion free.}

\bigskip\noindent
However, in order to state the conjecture on the center, we need some preliminaries.

\bigskip\noindent
The center of a group $G$ will be denoted by $Z(G)$. Let $\Gamma$ be a Coxeter graph, and let $\Gamma_1, \dots, \Gamma_l$ be the connected components of $\Gamma$. It is easily seen that 
\[
A_\Gamma = A_{\Gamma_1} \times \cdots \times A_{\Gamma_l}\,,
\]
hence 
\[
Z(A_\Gamma) = Z(A_{\Gamma_1}) \times \cdots \times Z(A_{\Gamma_l})\,.
\]
So, in order to understand the centers of the Artin-Tits groups, it suffices to consider connected Coxeter graphs. 

\bigskip\noindent
We say that a Coxeter graph $\Gamma$ is of {\it spherical type} if $W_\Gamma$ is finite. It is known that $Z(A_\Gamma) \simeq \Z$ if $\Gamma$ is connected and of spherical type (see \cite{Delig1, BriSai1}), and we believe that the center of any other Artin-Tits group associated to a connected Coxeter graph is trivial. In other words:

\bigskip\noindent
{\bf Conjecture B.}
{\it Let $\Gamma$ be a non-spherical connected Coxeter graph. Then $\Z(A_\Gamma)= \{1\}$.}

\bigskip\noindent
The third conjecture can be stated as easily as the first one:

\bigskip\noindent
{\bf Conjecture C.}
{\it Artin-Tits groups have solutions to the word problem.}

\bigskip\noindent
The statement of the last conjecture, known as {\it the $K(\pi,1)$ conjecture}, is more sophisticated than the previous ones. 

\bigskip\noindent
Let $\Gamma$ be a Coxeter graph.  The Coxeter group $W$ of $\Gamma$ has a faithful representation $W \hookrightarrow GL(V)$, called {\it canonical representation}, where $V$ is a real vector space of dimension $|S|$. The group $W$, viewed as a subgroup of $GL(V)$, is generated by reflections and acts properly discontinuously on a non-empty open cone $I$, called {\it Tits cone} (see \cite{Bourb1}). The set of reflections in $W$ is $R= \{w s w^{-1}; s \in S \text{ and } w \in W\}$, and $W$ acts freely and properly discontinuously on $I \setminus (\cup_{r \in R} H_r)$, where, for $r \in R$, $H_r$ denotes the hyperplane of $V$ fixed by $r$. Set
\[
M= M_\Gamma = (I \times V) \setminus \left( \bigcup_{r \in R} (H_r \times H_r) \right)\,.
\]
This is a connected manifold of dimension $2|S|$ on which the group $W$ acts freely and properly discontinuously. A key result in the domain is the following.

\bigskip\noindent
{\bf Theorem 2.1} 
(Van der Lek \cite{Lek1}). 
{\it The fundamental group of $M_\Gamma/W$ is isomorphic to the Artin-Tits group $A_\Gamma$.}

\bigskip\noindent
Recall that a CW-complex $X$ is a {\it a classifying space} for a (discrete) group $G$ if $\pi_1(X)=G$ and the universal cover of $X$ is contractible. These spaces play a prominent role in the calculation of the cohomology of groups (see \cite{Brown1}).

\bigskip\noindent
{\bf Conjecture D}
($K(\pi,1)$ conjecture).
{\it Let $\Gamma$ be a Coxeter graph. Then $M_\Gamma/W$ is a classifying space for $A_\Gamma$.}

\bigskip\noindent
Observe that, if Conjecture D holds for a Coxeter graph $\Gamma$, then the cohomological dimension of $A_\Gamma$ is finite, thus $A_\Gamma$ is torsion free (see \cite{Brown1}). In other words, Conjecture D implies Conjecture~A. 

\bigskip\noindent
Curiously, for each of the conjectures A, B, C, or D, the set of Artin-Tits groups for which a proof of the conjecture is known coincides (up to some exceptions) with the set of Artin-Tits groups for which proofs for the other three conjectures are known.

\bigskip\noindent
The first family of Artin-Tits groups which has been studied was the family of spherical type Artin-Tits groups (recall that $A_\Gamma$ is of spherical type if $W_\Gamma$ is finite). Conjectures A, B and C were proved for these groups in \cite{BriSai1,Delig1}, and Conjecture D was proved in \cite{Delig1}. 

\bigskip\noindent
We say that $\Gamma$ is of {\it large type} if $m_{s,t} \ge 3$ for all $s,t \in S$, $s\neq t$, and that $\Gamma$ is of {\it extra-large type} if $m_{s,t} \ge 4$ for all $s,t \in S$, $s \neq t$. Conjectures A and C were proved in \cite{AppSch1} for extra-large type Artin-Tits groups, and Conjecture B can be easily proved with the same techniques. Conjecture~D for large type Artin-Tits groups is a straightforward consequence of \cite{Hendr1}. 

\bigskip\noindent
Recall that a Coxeter graph $\Gamma$ is {\it free of infinity} if $m_{s,t} \neq \infty$ for all $s,t \in S$. We say that $\Gamma$ is of {\it FC type} if every free of infinity full subgraph of $\Gamma$ is of spherical type. On the other hand, we say that $\Gamma$ is of {\it dimension 2} if no full subgraph of $\Gamma$ with three or more vertices is of spherical type. 

\bigskip\noindent
After \cite{AppSch1} and \cite{Hendr1} the next significant step in the study of Artin-Tits groups and, more specifically, on the four conjectures stated above, was \cite{ChaDav1}. In this paper Conjecture D (and hence Conjecture~A) was proved for FC type Artin-Tits groups and 2-dimensional Artin-Tits groups (that include the Artin-Tits groups of large type). Later, Conjecture B was proved for Artin-Tits groups of FC type in \cite{Godel1}, and for 2-dimensional Artin-Tits groups in \cite{Godel2}.  On the other hand, Conjecture C was proved for Artin-Tits groups of FC type in \cite{Altob1, AltCha1}, and for 2-dimensional Artin-Tits groups in \cite{Cherma1}. Note also that a new solution to the word problem for FC type Artin-Tits groups will follow from Theorem~C. 

\bigskip\noindent
A challenging question in the domain is to prove Conjectures A, B, C, and D for the so-called Artin-Tits groups of affine type, that is, those Artin-Tits groups for which the associated Coxeter group is affine. As far as we know there is no known explicit algorithm for solving the word problem for these groups, except for the groups of type $\tilde A_n$ and $\tilde C_n$ (see \cite{Digne1, Digne2}). The techniques introduced in \cite{Allco1} may be used to solve the question for the groups of type $\tilde B_n$ and $\tilde D_n$, but we have no idea on how to  treat the groups of type $\tilde E_k$, $k=6,7,8$. The $K(\pi,1)$ conjecture was proved for the groups of type $\tilde A_n$ and $\tilde C_n$ in \cite{Okone1} (see also \cite{ChaPei1}) and for the groups of type $\tilde B_n$ in \cite{CaMoSa1}. The other cases are open. On the other hand, none of these papers addresses Conjecture~B. 


\section{From free of infinity Artin-Tits groups to Artin-Tits groups}

The aim of this section is to prove the following three theorems.

\bigskip\noindent
{\bf Theorem A.}
{\it If all free of infinity Artin-Tits groups are torsion free, then all Artin-Tits groups will be torsion free.}

\bigskip\noindent
{\bf Theorem B.}
{\it If, for every non-spherical connected free of infinity Coxeter graph $\Gamma$, we have $Z(A_\Gamma) = \{1\}$, then, for every non-spherical connected Coxeter graph $\Gamma$, we will have $Z(A_\Gamma) = \{1\}$.}

\bigskip\noindent
Note that, as pointed out before, if $\Gamma$ is a spherical type connected Coxeter graph, then $Z(A_\Gamma) \simeq \Z$ (see \cite{BriSai1}, \cite{Delig1}).

\bigskip\noindent
{\bf Theorem C.}
{\it If all free of infinity Artin-Tits groups have solutions to the word problem, then all Artin-Tits groups will have solutions to the word problem.}

\bigskip\noindent
In order to complete our discussion, we state the following theorem without proof. This is already proved in \cite{EllSko1} and \cite{GodPar1}.

\bigskip\noindent
{\bf Theorem D}
(Ellis, Sköldberg \cite{EllSko1}).
{\it If, for every free of infinity Coxeter graph $\Gamma$, the space $M_\Gamma/W_\Gamma$ is a classifying space for $A_\Gamma$, then, for every Coxeter graph $\Gamma$, the space $M_\Gamma/W_\Gamma$ will be a classifying space for $A_\Gamma$.}

\bigskip\noindent
The proofs of Theorems A, B and C use the notion of parabolic subgroup defined as follows. For $X\subset S$, we set $M_X = (m_{s,t})_{s,t \in X}$, we denote by $\Gamma_X$ the Coxeter graph of $M_X$, we denote by $W_X$ the subgroup of $W=W_\Gamma$ generated by $X$, we set $\Sigma_X = \{ \sigma_s; s \in X\}$, and we denote by $A_X$ the subgroup of $A = A_\Gamma$ generated by $\Sigma_X$. It is known that $(W_X,X)$ is the Coxeter system of $\Gamma_X$ (see \cite{Bourb1}), and that $(A_X,\Sigma_X)$ is the Artin-Tits system of $\Gamma_X$ (see \cite{Lek1, Paris1, GodPar1}). The subgroup $W_X$ is called {\it standard parabolic subgroup} of $W$, and $A_X$ is called {\it standard parabolic subgroup} of~$A$.

\bigskip\noindent
The proof of Theorem A is elementary. That of Theorem B is more complicated and requires some deeper knowledge on Artin-Tits groups. The proof of Theorem C is based on a combinatorial study of Artin-Tits groups made in \cite{GodPar1}. In all three cases we use the following observation.

\bigskip\noindent
{\bf Observation.}
Let $\Gamma$ be a Coxeter graph which is not free of infinity. Let $s,t \in S$, $s \neq t$, such that $m_{s,t} =\infty$. Set $X=S \setminus \{s\}$, $Y=S \setminus \{t\}$, and $Z = S \setminus \{s,t\}$. Recall that $(A_X, \Sigma_X)$ (resp. $(A_Y, \Sigma_Y)$, $(A_Z, \Sigma_Z)$) is an Artin-Tits system of type $\Gamma_X$ (resp. $\Gamma_Y$, $\Gamma_Z$), and that the natural homomorphisms $A_Z \to A_X$ and $A_Z \to A_Y$ are injective. Then the following equality is a direct consequence of the presentation of $A$.
\[
A = A_X \ast_{A_Z} A_Y\,.
\]

\bigskip\noindent
{\bf Proof of Theorem A.}
We assume that all free of infinity Artin-Tits groups are torsion free. Let $\Gamma$ be a Coxeter graph, $S$ its set of vertices, and $A=A_\Gamma$ its associated Artin-Tits group. We prove by induction on $|S|$ that $A$ is torsion free. If $|S|=1$, then $A$ is free of infinity, thus is torsion free. More generally, if $A$ is free of infinity then it is torsion free by hypothesis. Assume that $|S| \ge 2$ and $A$ is not free of infinity, plus the inductive hypothesis. Choose $s,t \in S$ such that $m_{s,t} = \infty$, and set $X=S \setminus \{s\}$, $Y=S \setminus\{t\}$ and $Z = S \setminus \{s,t\}$. Recall that $A=A_X \ast_{A_Z} A_Y$, hence, if $\alpha$ is a finite order element in $A$, then $\alpha$ is conjugate to an element in either $A_X$ or $A_Y$ (see \cite{Serre1}). By the inductive hypothesis $A_X$ and $A_Y$ are torsion free, thus $\alpha$ must be trivial.
\qed

\bigskip\noindent
The following two lemmas are preliminaries to the proof of Theorem B. 

\bigskip\noindent
{\bf Lemma 3.1.}
{\it Let $\Gamma$ be a connected spherical type Coxeter graph, let $S$ be its set of vertices, and let $X \subsetneq S$ be a proper subset of $S$. Then $Z(A) \cap A_X = \{1\}$.}

\bigskip\noindent
{\bf Proof.}
We assume first that $\Gamma$ is a spherical type Coxeter graph, but not necessarily connected. The {\it Artin-Tits monoid} of $\Gamma$ is defined by the monoid presentation
\[
A^+ = \langle \Sigma\ |\ \Pi(\sigma_s, \sigma_t:m_{s,t}) = \Pi (\sigma_t, \sigma_s : m_{s,t})\,,\ s,t \in S,\ s \neq t\text{ and } m_{s,t} \neq \infty \rangle^+\,.
\]
By \cite{BriSai1} the natural homomorphism $A^+ \to A$ is injective (see also \cite{Paris2}). Moreover, by \cite{Lek1}, we have $A^+ \cap A_X= A_X^+$ (this equality is also a direct consequence of the normal forms defined in \cite{Charn1}). On the other hand, it is easily deduced from the presentation of $A^+$ that, for $a \in A^+$, we have $a \in A_X^+$ if and only if any expression of $a$ in the elements of $\Sigma$ is actually a word in the elements of $\Sigma_X$. 

\bigskip\noindent
Now, we assume that $\Gamma$ is connected. Let $s_1, \dots, s_n$ be the vertices of $\Gamma$, and, for $i \in \{1, \dots, n\}$, set $\sigma_i = \sigma_{s_i}$. By \cite{BriSai1} the subgroup $Z(A)$ is the infinite cyclic group generated by
\[
\delta = (\sigma_1 \cdots \sigma_n)^h\,,
\]
where $h$ is the Coxeter number of $(W_\Gamma,S)$. Furthermore, $\delta$ does not depend on the choice of the ordering of $S$. Now, the previous observations clearly imply that $\delta^k \not\in A_X$ unless $k=0$, thus $Z(A) \cap A_X = \{1\}$.  
\qed

\bigskip\noindent
The next lemma is well-known and can be easily proved using normal forms in free products with amalgamation (see \cite{Serre1}).

\bigskip\noindent
{\bf Lemma 3.2.}
{\it Let $G=H_1 \ast_K H_2$ be the free product with amalgamation of two groups $H_1, H_2$ along a common subgroup $K$. Then $Z(G) = Z(H_1) \cap Z(H_2) \subset K$.}
\qed

\bigskip\noindent
{\bf Proof of Theorem B.}
We take a connected Coxeter graph $\Gamma$, and we prove Theorem B by induction on the number $|S|$ of vertices of $\Gamma$. If $|S|=1$ or, more generally, if $\Gamma$ is free of infinity, then, by \cite{BriSai1} and \cite{Delig1}, $Z(A_\Gamma) \simeq \Z$. 

\bigskip\noindent
Suppose $|S| \ge 2$ and $\Gamma$ is not free of infinity, plus the inductive hypothesis. Note that, in this case, $\Gamma$ is not of spherical type, thus we should prove that $Z(A)=\{1\}$. Let $s,t \in S$ be such that $m_{s,t} = \infty$. Set $X=S \setminus \{s\}$, $Y=S\setminus \{t\}$ and $Z = S \setminus \{s,t\}$. We have $A = A_X \ast_{A_Z} A_Y$ thus, by Lemma 3.2, 
\[
Z(A) = Z(A_X) \cap Z(A_Y) \subset A_Z\,.
\]

\bigskip\noindent
Let $X_1 \subset X$ be such that $\Gamma_{X_1}$ is the connected component of $\Gamma_X$ which contains $t$, let $X_2 = X \setminus X_1$, let $Y_1 \subset Y$ be such that $\Gamma_{Y_1}$ is the connected component of $\Gamma_Y$ which contains $s$, and let $Y_2 = Y \setminus Y_1$. First, we prove that $X_2 \subset Y_1$ and $Y_2 \subset X_1$.

\bigskip\noindent
Indeed, let $X_2' \subset X_2$ be such that $\Gamma_{X_2'}$ is a connected component of $\Gamma_{X_2}$. We have $X_2' \subset Y$, since $t \not \in X_2'$, and $\Gamma_{Y_1}$ is a connected component of $\Gamma_Y$, thus either $X_2' \subset Y_1$, or $X_2' \cap Y_1 = \emptyset$. Suppose $X_2' \cap Y_1 = \emptyset$, that is $X_2' \subset Y_2$. Let $x \in S \setminus X_2'$. If $x \neq s$, then $x \in X \setminus X_2'$, thus $m_{x,y}=2$ for all $y \in X_2'$, because $\Gamma_{X_2'}$ is a connected component of $\Gamma_X$. If $x=s$, then $m_{x,y}=2$ for all $y \in X_2'$, because $X_2' \subset Y_2$. This implies that $\Gamma_{X_2'}$ is a connected component of $\Gamma$: a contradiction. So, $X_2' \subset Y_1$. 

\bigskip\noindent
Set $Z_1 = X_1 \setminus \{t\}$. We have $A_Z = A_{Z_1} \times A_{X_2}$ and $A_{Z_1} \subset A_{X_1}$, thus
\[
Z(A_X) \cap A_Z = (Z(A_{X_1}) \cap A_{Z_1}) \times (Z(A_{X_2}) \cap A_{X_2}) = (Z(A_{X_1}) \cap A_{Z_1}) \times Z(A_{X_2})\,.
\]
If  $\Gamma_{X_1}$ is not of spherical type, then, by the inductive hypothesis, $Z(A_{X_1})=\{1\}$. If $\Gamma_{X_1}$ is of spherical type, then, by Lemma 3.2, $Z(A_{X_1}) \cap A_{Z_1} = \{1\}$. So, $Z(A_X) \cap A_Z = Z(A_{X_2}) \subset A_{X_2}$. Now, since $X_2 \subset Y_1$, we have 
\[
Z(A_Y) \cap A_{X_2} = (Z(A_{Y_1}) \cap A_{X_2}) \times (Z(A_{Y_2}) \cap \{1\}) = Z(A_{Y_1}) \cap A_{X_2}\,.
\]
If $\Gamma_{Y_1}$ is not of spherical type, then, by the inductive hypothesis, $Z(A_{Y_1})=\{1\}$. If $\Gamma_{Y_1}$ is of spherical type, then, since $X_2 \subsetneq Y_1$, we have $Z(A_{Y_1}) \cap A_{X_2} = \{1\}$ by Lemma 3.2. So, $Z(A_Y) \cap A_{X_2} = \{1\}$. Finally,
\[
Z(A) = Z(A_X) \cap Z(A_Y) = Z(A_X) \cap A_Z \cap Z(A_Y) \subset A_{X_2} \cap Z(A_Y) = \{1\}\,.
\qed
\] 

\bigskip\noindent
We turn now to the proof of Theorem C. we start with some definitions. Let $G = H_1 \ast_K H_2$ be an amalgamated free product. A {\it syllabic expression} of an element $\alpha \in G$ is a sequence $w=(\beta_1, \dots, \beta_l)$ of elements of $H_1 \cup H_2$ such that $\alpha = \beta_1 \beta_2 \cdots \beta_l$. We say that a syllabic expression $w'$ is an {\it elementary reduction of type I} of $w$ if there exists $i \in \{1, \dots, l\}$  such that $\beta_i=1$ and 
\[
w' = (\beta_1, \dots, \beta_{i-1},\beta_{i+1}, \dots, \beta_l)\,.
\]
We say that a syllabic expression $w'$ is an {\it elementary reduction of type II} of $w$ if there exists $i \in \{1, \dots, l-1\}$ and $j \in \{1,2\}$ such that $\beta_i,\beta_{i+1} \in H_j$ and
\[
w' = (\beta_1, \dots, \beta_{i-1},\beta_i \beta_{i+1}, \beta_{i+2}, \dots, \beta_l)\,.
\]
Implicit in this definition there is the assumption that every element of $K$ belongs to both, $H_1$ and $H_2$. Let $w,w'$ be two syllabic expressions of $\alpha$. If there is a finite sequence $w_0,w_1, \dots, w_n$ of syllabic expressions such that $w_0=w$, $w_n=w'$ and $w_i$ is an elementary reduction of $w_{i-1}$ for all $i \in \{1, \dots, n\}$, then we say that $w'$ is a {\it reduction} of $w$. A {\it reduced expression} is a syllabic expression which admits no elementary reduction. Clearly, a syllabic expression $w=(\beta_1, \dots, \beta_l)$ is reduced if and only if either $l=0$, or $l=1$ and $\beta_1 \neq 1$, or there exists a sequence $(k_1, \dots, k_l)$ in $\{ 1,2 \}$ such that $k_i \neq k_{i+1}$ for all $i \in \{1, \dots, l-1\}$, and $\beta_i \in H_{k_i} \setminus K$ for all $i \in \{1, \dots, l\}$. The following theorem is well-known and widely used in the study of amalgamated free products of groups (see \cite{Serre1}).

\bigskip\noindent
{\bf Theorem 3.3.}
{\it Let $G = H_1 \ast_K H_2$ be an amalgamated free product. Let $\alpha \in G$, and let $w=(\beta_1, \dots, \beta_l)$, $w'=(\beta_1', \dots, \beta_k')$ be two reduced syllabic expressions of $\alpha$. Then $k=l$ and there exist $\gamma_0,\gamma_1, \dots, \gamma_l \in K$ such that
\[
\gamma_0=\gamma_l=1,\quad \text{and } \beta_i'\gamma_i= \gamma_{i-1} \beta_i \text{ for all } i \in \{1, \dots, l\}\,.
\]
In particular, $\alpha=1$ if and only if its (unique) reduced expression is the empty sequence of length $0$.
\qed}

\bigskip\noindent
For a given set $X$ we denote by $X^\ast$ the free monoid generated by $X$. If $G$ is a group generated by a family $S$ and $w \in (S \sqcup S^{-1})^\ast$, we denote by $[w]_S$ the element of $G$ represented by $w$. On the other hand, we will assume for the next corollary that $H_1,H_2$, and $K$ are finitely generated, and we take finite generating sets $S_1,S_2$, and $T$, for $H_1,H_2$, and $K$, respectively. 

\bigskip\noindent
{\bf Corollary 3.4.}
{\it Assume that, for each $j \in \{1,2\}$,
\begin{itemize}
\item[{\bf (a)}]
$H_j$ has a solution to the word problem; and
\item[{\bf (b)}]
there exists an algorithm which, given a word $w \in (S_j \sqcup S_j^{-1})^\ast$, decides whether $\beta=[w]_{S_j}$ belongs to $K$ and, if yes, determines a word $u \in (T \sqcup T^{-1})^\ast$ which represents $\beta$.
\end{itemize}
Then $G=H_1 \ast_K H_2$ has a solution to the word problem.}

\bigskip\noindent
{\bf Proof.}
Let $\alpha \in G$. A {\it syllabic pre-expression} of $\alpha$ is defined to be a pair of sequences of equal length, 
\[
U=((w_1, \dots, w_l),(j_1, \dots, j_l))\,,
\]
such that $j_i \in \{1,2\}$, and $w_i \in (S_{j_i} \sqcup S_{j_i}^{-1})^\ast$, for all $i \in \{1, \dots, l\}$, and $\alpha=[w_1 \cdots w_l]_{S_1 \cup S_2}$. The {\it syllabic realization} of $U$ is the syllabic expression
\[
[U] = (\beta_1, \dots, \beta_l)\,,
\]
where $\beta_i = [w_i]_{S_{j_i}}$ for all $i \in \{1, \dots, l\}$. We will say that $U$ is {\it reduced} if its realization $[U]$ is reduced. 

\bigskip\noindent
Let $U=((w_1, \dots, w_l),(j_1, \dots, j_l))$ be a syllabic pre-expression of $\alpha$. Suppose there exists $i \in \{1, \dots, l\}$ such that $[w_i]_{S_{j_i}}=1$. Note that, Since $H_{j_i}$ has a solution to the word problem, there is an algorithm which decides whether this is true. Set
\[
U' = ((w_1, \dots, w_{i-1},w_{i+1}, \dots, w_l), (j_1, \dots, j_{i-1},j_{i+1}, \dots, j_l))\,.
\]
Then $U'$ is a syllabic pre-expression of $\alpha$ and $[U']$ is an elementary reduction of type I of $[U]$. Suppose there exists $i \in \{1, \dots, l-1\}$ such that $j_i = j_{i+1}$. Set
\[
U' = ((w_1, \dots, w_{i-1},w_iw_{i+1},w_{i+2}, \dots, w_l), (j_1, \dots, j_{i-1},j_i, j_{i+2}, \dots, j_l))\,.
\]
Then $U'$ is a syllabic pre-expression of $\alpha$ and $[U']$ is an elementary reduction of type II of $[U]$. Suppose there exists $i \in \{1, \dots, l-1\}$ such that $j_i \neq j_{i+1}$, but $\beta_i = [w_i]_{S_{j_i}} \in K$. Recall that, by hypothesis, there is an algorithm which decides whether this is true. Again by hypothesis, there is an algorithm which determines a word $u_i \in (T \sqcup T^{-1})^\ast$ which represents $\beta_i$. Writing the elements of $T$ in the generators $S_{j_{i+1}} \sqcup S_{j_{i+1}}^{-1}$ of $H_{j_{i+1}}$, we get from $u_i$ a word $w_i' \in (S_{j_{i+1}} \sqcup S_{j_{i+1}}^{-1})^\ast$ which represents $\beta_i$. Set
\[
U' = ((w_1, \dots, w_{i-1}, w_i'w_{i+1}, w_{i+2}, \dots, w_l),(j_1, \dots, j_{i-1},j_{i+1},j_{i+2}, \dots, j_l))\,.
\]
Then $U'$ is a syllabic pre-expression of $\alpha$ and $[U']$ is an elementary reduction of type II of $[U]$. Suppose there exists $i \in \{1, \dots, l-1\}$ such that $j_i \neq j_{i+1}$, but $\beta_{i+1} = [w_{i+1}]_{S_{j_{i+1}}} \in K$. Then, in the same manner as in the previous case, one can explicitely calculate an expression $w_{i+1}'$ for $\beta_{i+1}$ in $(S_{j_i} \sqcup S_{j_i}^{-1})^\ast$,
\[
U' = ((w_1, \dots, w_{i-1}, w_iw_{i+1}', w_{i+2}, \dots, w_l),(j_1, \dots, j_{i-1},j_i, j_{i+2}, \dots, j_l))\,
\]
is a syllabic pre-expression of $\alpha$, and $[U']$ is an elementary reduction of type II of $[U]$. It is easily seen that all elementary reductions of $[U]$ are of these forms. Hence, from a given syllabic pre-expression $U$ of $\alpha$ one can effectively calculate a reduced syllabic pre-expression $U'$ of $\alpha$. Moreover, by Theorem 3.3, we have $\alpha=1$ if and only if $U'$ is a pair of empty sequences (of length $0$).

\bigskip\noindent
Let $w=s_1^{\varepsilon_1} \cdots s_l^{\varepsilon_l} \in (S_1 \sqcup S_2 \sqcup S_1^{-1} \sqcup S_2^{-1})^\ast$, and let $\alpha=[w]_{S_1 \cup S_2} \in G$. For $i \in \{1, \dots, l\}$ we set $j_i=1$ if $s_i \in S_1$ and $j_i=2$ if $s_i \in S_2$. Then
\[
U=((s_1^{\varepsilon_1}, \dots, s_l^{\varepsilon_l}),(j_1, \dots, j_l))
\]
is a syllabic pre-expression of $\alpha$, and, by the above, one can decide from $U$ if $\alpha=1$.
\qed

\bigskip\noindent
Now, the key in the proof of Theorem C is the following result proved in \cite{GodPar1}.

\bigskip\noindent
{\bf Theorem 3.5}
(Godelle, Paris \cite{GodPar1}).
{\it Let $\Gamma$ be a Coxeter graph, let $S$ be be its set of vertices, let $X \subset S$, and let $(A,\Sigma)$ be the Artin-Tits system of $\Gamma$. If $A$ has a solution to the word problem, then there is an algorithm which, given a word $w \in (\Sigma \sqcup \Sigma^{-1})^\ast$, decides whether $\alpha = [w]_\Sigma$ belongs to $A_X$, and, if yes, determines a word $w' \in (\Sigma_X \sqcup \Sigma_X^{-1})^\ast$ which represents $\alpha$.}
\qed

\bigskip\noindent
{\bf Proof of Theorem C.}
We assume that all free of infinity Artin-Tits groups have solutions to the word problem. Let $\Gamma$ be a Coxeter graph, let $S$ be its set of vertices, and let $(A,\Sigma)$ be its associated Artin-Tits system. We prove by induction on $|S|$ that $A$ has a solution to the word problem.

\bigskip\noindent
If $|S|=1$ or, more generally, if $\Gamma$ is free of infinity, then, by hypothesis, $A$ has a solution to the word problem. Assume that $\Gamma$ is not free of infinity, plus the inductive hypothesis. Let $s,t \in S$, $s \neq t$, such that $m_{s,t} = \infty$. Set $X=S \setminus \{s\}$, $Y=S \setminus \{t\}$, and $Z=S \setminus \{s,t\}$. We have $A=A_X \ast_{A_Z} A_Y$ and, by the inductive hypothesis, $A_X,A_Y$ have solutions to the word problem. Moreover, by Theorem 3.5, for $T=X$ or $Y$, there is an algorithm which, given a word $w \in (\Sigma_T \sqcup \Sigma_T^{-1})^\ast$, decides whether $\alpha=[w]_{\Sigma_T}$ lies in $A_Z$ and, if yes, determines a word $w' \in (\Sigma_Z \sqcup \Sigma_Z^{-1})^\ast$ which represents $\alpha$. We conclude by Corollary 3.4 that $A$ has a solution to the word problem.
\qed



\bigskip\bigskip\noindent
{\bf Eddy  Godelle,}

\smallskip\noindent
Université de Caen, Laboratoire LMNO, UMR 6139 du CNRS, Campus II, 14032 Caen cedex, France. 

\smallskip\noindent
E-mail: {\tt eddy.godelle@unicaen.fr}

\bigskip\noindent
{\bf Luis Paris,}

\smallskip\noindent 
Université de Bourgogne, Institut de Mathématiques de Bourgogne, UMR 5584 du CNRS, B.P. 47870, 21078 Dijon cedex, France.

\smallskip\noindent
E-mail: {\tt lparis@u-bourgogne.fr}



\begin{thebibliography}{99}

\small

\bibitem{Allco1}
{\bf D. Allcock.}
{\it Braid pictures for Artin groups.}
Trans. Amer. Math. Soc. {\bf 354} (2002), no. 9, 3455--3474. 

\bibitem{Altob1}
{\bf J. Altobelli.}
{\it The word problem for Artin groups of FC type.}
J. Pure Appl. Algebra {\bf 129} (1998), no. 1, 1--22.

\bibitem{AltCha1}
{\bf J. Altobelli, R. Charney.}
{\it A geometric rational form for Artin groups of FC type.}
Geom. Dedicata {\bf 79} (2000), no. 3, 277--289.

\bibitem{AppSch1}
{\bf K.\,I. Appel, P.\,E. Schupp.}
{\it Artin groups and infinite Coxeter groups.}
Invent. Math. {\bf 72} (1983), no. 2, 201--220.

\bibitem{Bourb1}
{\bf N. Bourbaki.}
{\it Eléments de mathématique. Fasc. XXXIV. Groupes et algèbres de Lie. Chapitre IV: Groupes de Coxeter et systèmes de Tits. Chapitre V: Groupes engendrés par des réflexions. Chapitre VI: Systèmes de racines.}
Actualités Scientifiques et Industrielles, No. 1337, Hermann, Paris, 1968. 

\bibitem{BriSai1}
{\bf E. Brieskorn, K. Saito.}
{\it Artin-Gruppen und Coxeter-Gruppen.}
Invent. Math. {\bf 17} (1972), 245--271. 

\bibitem{Brown1}
{\bf K.\,S. Brown.}
{\it Buildings.}
Springer-Verlag, New York, 1989.

\bibitem{CaMoSa1}
{\bf F. Callegaro, D. Moroni, M. Salvetti.}
{\it The $K(\pi,1)$ problem for the affine Artin group of type $\tilde B_n$ and its cohomology.}
J. Eur. Math. Soc. (JEMS) {\bf 12} (2010), no. 1, 1--22.

\bibitem{Charn1}
{\bf R. Charney.}
{\it Geodesic automation and growth functions for Artin groups of finite type.}
Math. Ann. {\bf 301} (1995), no. 2, 307--324.

\bibitem{ChaDav1}
{\bf R. Charney, M.\,W. Davis.}
{\it The $K(\pi,1)$-problem for hyperplane complements associated to infinite reflection groups.}
J. Amer. Math. Soc. {\bf 8} (1995), no. 3, 597--627.

\bibitem{ChaPei1}
{\bf R. Charney, D. Peifer.}
{\it The $K(\pi,1)$-conjecture for the affine braid groups.}
Comment. Math. Helv. {\bf 78} (2003), no. 3, 584--600.

\bibitem{Cherma1}
{\bf A. Chermak.}
{\it Locally non-spherical Artin groups.}
J. Algebra {\bf 200} (1998), no. 1, 56--98.

\bibitem{Davis1}
{\bf M.\,W. Davis.}
{\it The geometry and topology of Coxeter groups.}
London Mathematical Society Monographs Series, 32. Princeton University Press, Princeton, NJ, 2008.

\bibitem{Delig1}
{\bf P. Deligne.}
{\it Les immeubles des groupes de tresses généralisés.}
Invent. Math. {\bf 17} (1972), 273--302.

\bibitem{Digne1} 
{\bf F. Digne.}
{\it Présentations duales des groupes de tresses de type affine $\widetilde A$.}
Comment. Math. Helv. {\bf 81} (2006), no. 1, 23--47.

\bibitem{Digne2}
{\bf F. Digne.}
{\it A Garside presentation for Artin-Tits groups of type $\tilde C_n$.}
Preprint, arXiv: 1002.4320. 

\bibitem{EllSko1}
{\bf G. Ellis, E. Sköldberg.}
{\it The $K(\pi,1)$ conjecture for a class of Artin groups.}
Comment. Math. Helv. {\bf 85} (2010), no. 2, 409--415. 

\bibitem{Godel1}
{\bf E. Godelle.}
{\it Parabolic subgroups of Artin groups of type FC.}
Pacific J. Math. {\bf 208} (2003), no. 2, 243--254. 

\bibitem{Godel2}
{\bf E. Godelle.}
{\it Artin-Tits groups with CAT(0) Deligne complex.}
J. Pure Appl. Algebra {\bf 208} (2007), no. 1, 39--52.

\bibitem{GodPar1}
{\bf E. Godelle, L. Paris.}
{\it $K(\pi,1)$ and word problems for infinite type Artin-Tits groups, and applications to virtual braid groups.}
Preprint, arXiv:1007.1365. 

\bibitem{Hendr1}
{\bf H. Hendriks.}
{\it Hyperplane complements of large type.}
Invent. Math. {\bf 79} (1985), no. 2, 375--381. 

\bibitem{Humph1}
{\bf J.\,E. Humphreys.}
{\it Reflection groups and Coxeter groups.}
Cambridge Studies in Advanced Mathematics, 29. Cambridge University Press, Cambridge, 1990. 

\bibitem{Lek1}
{\bf H. van der Lek.}
{\it The homotopy type of complex hyperplane complements.}
Ph. D. Thesis, Nijmegen, 1983.

\bibitem{Okone1}
{\bf C. Okonek.}
{\it Das $K(\pi ,\,1)$-Problem für die affinen Wurzelsysteme vom Typ $A_{n}$, $C_{n}$.}
Math. Z. {\bf 168} (1979), no. 2, 143--148. 

\bibitem{Paris1}
{\bf L. Paris.}
{\it Universal cover of Salvetti's complex and topology of simplicial arrangements of hyperplanes.}
Trans. Amer. Math. Soc. {\bf 340} (1993), no. 1, 149--178.

\bibitem{Paris2}
{\bf L. Paris.}
{\it Artin monoids inject in their groups.}
Comment. Math. Helv. {\bf 77} (2002), no. 3, 609--637. 

\bibitem{Serre1}
{\bf J.-P. Serre.}
{\it Arbres, amalgames, ${\rm SL}_{2}$.} 
Astérisque, No. 46. Société Mathématique de France, Paris, 1977. 

\bibitem{Tits1}
{\bf J. Tits.}
{\it Groupes et géométries de Coxeter.}
Institut des Hautes Etudes Scientifiques, Paris, 1961.

\bibitem{Tits2}
{\bf J. Tits.}
{\it Normalisateurs de tores. I. Groupes de Coxeter étendus.}
J. Algebra {\bf 4} (1966), 96--116.

\end{thebibliography}
\end{document}